\documentclass[11pt]{amsart}
\usepackage{amsmath,amsthm,amsfonts,amssymb,bm,graphicx}
\usepackage
    {hyperref,cleveref}
\hypersetup{colorlinks=true,citecolor=blue,linkcolor=blue,urlcolor=blue,
pdfstartview=FitH }

\theoremstyle{plain}
\newtheorem{theorem}{Theorem}[section]

\theoremstyle{remark}
\newtheorem{rem}[theorem]{Remark}
\newtheorem*{acknowledgements}{Acknowledgements}

\numberwithin{equation}{section}

\numberwithin{table}{section}

\newcommand{\qq}{\mathfrak{q}}
\newcommand{\C}{\mathbb{C}}
\newcommand{\Fq}{\mathbb{F}_q}
\newcommand{\R}{\mathbb{R}}

\newcommand{\Z}{\mathbb{Z}}
\newcommand{\N}{\mathbb{N}}
\newcommand{\Gri}{\ensuremath{\mathcal{O}}}
\newcommand{\GL}{\mathsf{GL}}
\newcommand{\SL}{\mathsf{SL}}

\begin{document}

\author{Valentin Blomer}
\address{Mathematisches Institut, Endenicher Allee 60, 53115 Bonn, Germany} \email{blomer@math.uni-bonn.de}
 
\author{Christopher Voll}
\address{Fakult\"at f\"ur Mathematik, Universit\"atsstra{\ss}e 25
33615 Bielefeld, Germany}\email{voll@math.uni-bielefeld.de}

 \title{Analytic properties of representation zeta functions of groups
   of type $\mathsf{A_2}$}



\begin{abstract}
 We study analytic properties of the representation zeta functions of
 arithmetic groups of type $\mathsf{A}_2$, such as $\SL_3(\Z)$. In
 particular, we uncover further poles of these functions and determine
 a natural boundary for their meromorphic continuation beyond their
 abscissa of convergence. We analyse both the number field and
 function field case. 
\end{abstract}


\setcounter{tocdepth}{2}  \maketitle 

\maketitle


\thispagestyle{empty}

\section{Introduction}

\subsection{Euler products: meromorphic continuations and
  natural boundaries} $L$-functions in the Selberg class have an Euler
product and an analytic continuation to the entire complex plane with a
pole at most at $s = 1$.  Many natural $L$-functions come as an Euler
product because of some underlying local-global principle, but the
analytic continuation is by no means clear, either because it is hard
to prove (which is one of the major obstacles in the Langlands
program), or because it is simply not true.

There is no reason that a random Euler product should have analytic
continuation, and even if it has some structure, this defines rarely
an entire function. The prototypical result goes back to Estermann
\cite{Es} who proved that an Euler product of the shape $\prod_p
h(p^{-s})$ with $h \in \Bbb{Z}[x]$ satisfying $h(0) = 1$ has a natural
boundary at $\Re s = 0$ unless $h$ is a product of cyclotomic
polynomials. This has been generalized in various ways, most notably
in \cite[Ch.~5]{dSW} to Euler products of the form
\begin{equation}\label{poly}
    \prod_p h(p, p^{-s})
 \end{equation}   
for a bivariate polynomial~$h\in \Bbb{Z}[x, y]$. While the theory
developed in \cite{dSW} is not exhaustive, it indicates that such
Euler products typically have---provably or conjecturally---a natural
boundary a little bit to the left of the abscissa of absolute
convergence, and there should be a recipe to read it off from the
shape of~$h$. Nevertheless, the existence and location of the natural
boundary is rather subtle, cf.\ e.g.\ the discussion in \cite{BSP},
and remains unknown in many seemingly simple cases. See also \cite{Al} for a survey. 

\subsection{Euler products from representation growth of groups}\label{subsec:EPfromRG}
In this note we consider analytic properties of Euler products yet
more complicated than \eqref{poly}, arising as representation zeta
functions of arithmetic groups. For a group $G$ and $n \in \Bbb{N}$
let $r_n(G)$ denote the number of inequivalent $n$-dimensional
irreducible complex (continuous, if $G$ is topological)
representations of~$G$. We call $G$ (\emph{representation})
\emph{rigid} if $r_n(G)$ is finite for all~$n$. In this case, the
\emph{representation zeta function} $$\zeta_G(s) = \sum_{n=1}^\infty
r_n(G) n^{-s}$$ of $G$ converges for some $s \in \Bbb{C}$ if and only
if $G$ has \emph{polynomial representation growth}, i.e.\ $r_n(G)$ is
bounded by a polynomial in~$n$. This holds for certain arithmetic
groups of type ${\sf A}_2$ defined over number fields, such
as~$\SL_3(\Z)$. The  zeta functions of these groups were intensely studied
in~\cite{AKOV}.

We proceed to describe the precise set-up, recalling well-known facts,
e.g.\ from~\cite[\S~1.1]{AKOV}. Throughout this paper, let $k$ be a
number field (i.e.\ a finite extension of $\Bbb{Q}$) or a function field (i.e.\ a finite extension of $\Bbb{F}_{\qq}(T)$) with ring of integers~$\Gri$. For a place $v$ of $k$, we
write $k_v$ for the completion of $k$ at~$v$ and, if $v$ is
non-Archimedean, $\Gri_v$ for the completion of $\Gri$ at~$v$.  Let
$S$ be a finite set of places of $k$, including all the Archimedean
ones in the number field case, and let $\Gri_S = \{x \in k \mid   x \in
\Gri_v \text{ for all } v\not\in S\}$ be the ring of $S$-integers in~$k$.

Let $\mathbf{H}$ be a connected, simply-connected absolutely almost
simple algebraic group defined over~$k$, with a fixed embedding
into~$\GL_d$ for some~$d \in \N$. We consider the arithmetic group
$\mathbf{H}(\Gri_S) = \mathbf{H}(k) \cap \GL_d(\Gri_S)$.  If
$\mathbf{H}(\Gri_S)$ has the strong Congruence Subgroup Property
(i.e.\ the congruence kernel ${\rm ker}(\widehat{\mathbf{H}(\Gri_S)}
\rightarrow \mathbf{H}(\widehat{\Gri_S}))$ is trivial), abbreviated by
sCSP, then we have
\begin{equation}\label{equ:EP}
  \zeta_{\mathbf{H}(\Gri_S)}(s) = \zeta_{\mathbf{H}(\C)}(s)^{r}
  \prod_{v\not\in S} \zeta_{\mathbf{H}(\Gri_v)}(s),
\end{equation}
where $r$ is the degree of $k$ over $\Bbb{Q}$ in the number field case and $0$ in the function field case; 
cf.\ \cite[Proposition~1.3]{LL}. 
The Archimedean factors $\zeta_{\mathbf{H}(\C)}$
enumerate the finite-dimensional, irreducible rational representations
of the algebraic group $\mathbf{H}(\C)$; their contribution to the
Euler product reflects Margulis super-rigidity. The non-Archimedean
Euler factors indexed by the places not in $S$ are all rational
functions, albeit not just in $q_v$ and $q_v^{-s}$, where $q_v$
denotes the residue field cardinality at~$v$. Computing these rational
functions has proven to be very challenging. Explicit formulae seem
only to be known for groups of (Lie) type $\mathsf{A}_2$
and~$\mathsf{A}_1$.

\subsection{Groups of type $\mathsf{A}_2$}
Assume now that $\mathbf{H}$ is a connected, simply-connected,
absolutely almost simple algebraic group of type ${\sf A}_2$ defined
over a number field $k$ or over a function field $k$ with
characteristic greater than 3.  We consider the situation when it is
either an inner form arising from a matrix algebra over a central
division algebra over $k$ or an outer form over a central division
algebra over a quadratic extension $K/k$ (with the same field
$\Bbb{F}_{\mathfrak{q}}$ of constants in the function field case so
that $K/k$ is a ``geometric'' extension).

Assume that $\mathbf{H}(\Gri_S)$ has the sCSP, so \eqref{equ:EP}
applies. It is known that the abscissa of convergence of
$\zeta_{\mathbf{H}({\Gri_S})}$ is equal to $1$ and that it has
meromorphic continuation to $\Re s > 1-\delta$ for some $\delta > 0$
(specifically, $\delta = 1/6$ if $k$ is a number field) with a double
pole at $s=1$ and no pole in $1 - \delta < \Re s < 1$; see
\cite[Thm.~A]{OPS/25}.

 In this paper we prove the following best possible refinement of
 these results, uniformly for number fields and function fields.  Let
 $\zeta_3 \in \overline{k}$ denote a primitive third root of unity. If
 $k$ is a number field we distinguish two cases:
$$\begin{cases} \text{\emph{Case (A)}:} & \text{($\mathbf{H}$ is an
     outer form and $K = k(\zeta_3)$)} \,\text{ or }\,
   \text{($\mathbf{H}$ is an inner form and $\zeta_3 \in
     k$)},\\ \text{\emph{Case (B)}:} & \text{otherwise.}\end{cases}$$
 For the investigation of $\zeta_{{\rm SL}_3(\Bbb{Z})}$, for instance,
 Case~(B) applies, since we are in the case of an inner form and
 $\zeta_3 \not \in \Bbb{Q}$.

If $k$ is a function field with field of constants $\Bbb{F}_{\mathfrak{q}}$, we define $\eta = 0$ if $\zeta_3 \in k$, equivalently $\mathfrak{q} \equiv 1$ (mod 3), and $\eta = 1$ if $\zeta_3 \not\in k$, equivalently $\mathfrak{q} \equiv -1$ (mod 3).  Here we distinguish the cases (A): $\mathbf{H}$ is an inner form and (B): $\mathbf{H}$ is an outer form. 


\begin{theorem}\label{thm1}
  Assume that $\mathbf{H}(\Gri_S)$ is an arithmetic group of type
  $\mathsf{A}_2$ defined as above over a global field $k$ with characteristic 0 or greater than 3, as described above, satisfying the sCSP. Then the
  function 
  $\zeta_{\mathbf{H}({\Gri_S})}$ has meromorphic continuation to
  $\Re s > 5/8$ and a natural boundary at $\Re s = 5/8$.
  
  If $k$ is a number field, it has a
  double pole at $s=1$ and a pole at $s = 4/5$ of order 9 in Case (A)  and of order 5 in Case (B), and it has no other poles in $\Re s > 3/4$. 
  
  If $k$ is a function field with field of constants $\Bbb{F}_{\mathfrak{q}}$, it has    double poles at $s \in 1 + \frac{2\pi i}{\log \qq} \Bbb{Z}$ and poles at $s \in 4/5 + \frac{\pi i}{\log \qq}  (\eta + 2\Bbb{Z}) $ of order 9 in Case (A)  and of order 5 in Case (B), and it has no other poles in $\Re s > 3/4$.

\end{theorem}

\begin{rem}\
  \begin{enumerate}
    \item 
 That the degree of representation growth, i.e.\ the abscissa of
 convergence of $\zeta_{\mathbf{H}(\Gri_S)}$, only depends on the Lie
 type of the abstract group $\mathbf{H}$ (and not, for instance, on
 the ring $\Gri_S$) is an instance of a more general phenomenon: the
 degree of representation growth of an arithmetic group of the form
 $\mathbf{H}(\Gri_S)$ with the (weak) Congruence Subgroup Property
 (i.e.\ with finite congruence kernel) only depends on the root system
 associated with the algebraic group $\mathbf{H}$; see
 \cite[Thm.~1.1]{AKOV4}. As similar invariance phenomenon is
 established in \cite{DV/17} for representation zeta functions
 associated with unipotent group schemes.

 In contrast, we see here for the first time a situation where finer
 invariants, such as the pole order of the second right-most pole of
 $\zeta_{\mathbf{H}(\Gri_S)}$, here at $s=4/5$, depend more subtly on
 the underlying arithmetic structures.

  \item It would be interesting to have a conceptual explanation for
    the constant $5/8$. The local Euler factors come naturally as
    finite sums, and a trivial analysis of each summand would only
    lead to analytic continuation up to $\Re s > 2/3$. However, there
    is substantial cancellation---not algebraically, but
    asymptotically---which allows to continue to $\Re s > 5/8$, but
    not further; cf.\ also Remark \ref{rem21}.
    \end{enumerate}
\end{rem}


The proof of \Cref{thm1} is given in~\Cref{subsec:proof.A2}. It is
based on an analysis of the non-Archimedean factors in \eqref{equ:EP}
provided in~\cite{AKOV}. The Witten zeta function
$\zeta_{\mathbf{H}(\C)}$ poses no particular difficulty, so we turn
directly to the non-Archimedean part of this Euler product. As noted
above, all its factors are rational functions. We recall the explicit
yet intricate formulae for almost all of these functions
in~\Cref{subsubsec:non-arch}. The product over their denominators
poses no difficulty, as it yields a product of two translates of the
zeta function $\zeta_k$ of~$k$, with a few Euler factors omitted. The
challenge thus lies in the analysis of the product of the
numerators. This Euler product is not of the form~\eqref{poly}, for
which the theory of \cite{dSW} would be available \emph{mutatis
mutandis}, i.e.\ with the product running over places $v\not\in S$ and
$p$ replaced by the residue field cardinality~$q$. Instead, we are led
to consider Euler products of polynomials in $q^{-s}$ and $n_i^{-s}$,
for polynomial expressions $n_i$ in $q$, with coefficients depending
(mildly) on invariants of the arithmetic group $\mathbf{H}$ and
certain congruence classes of the residue field cardinality~$q$.  The
remaining ``exceptional'' Euler factors have little bearing on this
analysis; we discuss them, together with the Archimedean factors, in
\Cref{subsubsec:except}. The hardest part of the proof is the analysis
of the natural boundary. At this point the argument in the function
field case diverges from the number field case.

We are able to detect a second order term in the asymptotic formula
for the Dirichlet series coefficients $r_n(\mathbf{H}({\Gri_S}))$
of~$\zeta_{\mathbf{H}({\Gri}_S)}$. In the number field case we need to
insert a smooth weight for reasons to be explained in a moment.

\begin{theorem}\label{thm2}
  Assume that $\mathbf{H}(\Gri_S)$ is as in~\Cref{thm1}. There exist
  polynomials $P, \tilde{P}\in\R[X]$ with $\deg P = 1$ and $\deg
  \tilde{P} = 8$ in Case (A) and $\deg \tilde{P} = 4$ in Case (B), so
  that for every $\epsilon > 0$ the following hold:
  
  If $k$ is a number field, then
  \begin{equation}\label{sim:smooth}
    \sum_{n=1}^\infty r_n(\mathbf{H}({\Gri_S})) e^{-n/x} =x P(\log x)
    + x^{4/5} \tilde{P}(\log x) + O_{\epsilon}(x^{3/4 + \epsilon}).
  \end{equation}
  
  If $k$ is a function field with constant field $\Bbb{F}_{\qq}$, then 
   \begin{equation*}
     r_{\qq^n}(\mathbf{H}({\Gri_S}))  =\qq^n P(\log \qq^n)
    + \qq^{\frac{4}{5}n} \tilde{P}(\log \qq^n) + O_{\epsilon}(\qq^{n(3/4 + \epsilon)}).
  \end{equation*}

\end{theorem}

\begin{rem}
That we cannot expect a power-saving in a sharp cut-off count in the situation of \eqref{sim:smooth} is not
surprising, as an inspection of the Euler product over the factors
\eqref{euler} shows: it is at least as difficult as the basic function
$\sum_n \phi(n)^{-s} = \zeta(s) \prod_p ( 1 + (p-1)^{-s} - p^{-s}) $
where at the current state of knowledge for the summatory function $\sum_{\phi(n) \leq x} 1$ no
power-saving is available; see~\cite{Ba}.
\end{rem}

\begin{rem}
  The identity \eqref{sim:smooth} invites a comparison with the
  asymptotic statement that
  \begin{equation}\label{equ:ass}
    \sum_{n=1}^x r_n(\mathbf{H}({\Gri_S})) \sim
    c(\mathbf{H}(\Gri_S))\cdot x \log x \quad \textup{ for }x
    \rightarrow \infty
  \end{equation} for a constant
    $c(\mathbf{H}(\Gri_S))\in\R_{>0}$; see
  \cite[Cor.~B(2)]{AKOV}. Comparing \eqref{equ:ass} with
  \eqref{sim:smooth} yields that $c(\mathbf{H}(\Gri_S))$ is the leading
  coefficient of $P_1$. The constant $c(\SL_3(\Z))$ is discussed
  in~\cite[Sec.~7.1]{AKOV}.
   \end{rem}

\subsection{Groups of type $\mathsf{A}_1$}
Assume now that $\mathbf{H}$ is of type $\mathsf{A}_1$, viz.\ a form
of $\SL_2$, and that $S$ contains all places dividing~$2$
and~$\infty$. The group $\SL_2(\Z)$ does not have the strong Congruence
Subgroup Property, but groups of the form $\SL_2(\Gri_S)$ do, for
sufficiently large finite sets of places~$S$. Zeta functions of groups
of the form $\mathbf{H}({\Gri_S})$ with the sCSP have been considered
e.g.\ in \cite[\S~10]{LL}. By \cite[Thm.~10]{LL} we know that
$\zeta_{\mathbf{H}(\Gri_S)}$ has abscissa of convergence equal to
$2$. A porism of their result is that the Euler product \eqref{equ:EP}
allows for some meromorphic continuation, unveiling a simple pole at
$s=2$. We extend these results as follows.

\begin{theorem}\label{prop3}
  Assume that $\mathbf{H}(\Gri_S)$ is an arithmetic group of type
  $\mathsf{A}_1$ defined as above, over a number field or over a function field of characteristic greater than 3, satisfying the sCSP.  The function
  $\zeta_{\mathbf{H}({\Gri_S})}$ has meromorphic continuation to
  $\Re s > 1$ with a simple pole at $s=2$ (resp.\ simple poles at $s = 2 + \frac{2\pi i}{\log \qq} n$, $n \in\Bbb{Z}$, if $k$ is  a function field with constant field $\Bbb{F}_{\qq}$) and no further poles. It has a branch cut singularity at $s=1$. 

There exists a constant $c > 0$ such that for every $\epsilon>0$ the
following hold:
  
  If $k$ is a number field, then
$$\sum_{n=1}^\infty r_n(\mathbf{H}({\Gri_S})) e^{-n/x} = c x^2   + O_{\epsilon}(x^{1+\epsilon}).$$

If $k$ is a function field with constant field $\Bbb{F}_{\qq}$, then
$$  r_{\qq^n} (\mathbf{H}({\Gri_S}))  = c \qq^{2n}   + O_{\epsilon}(\qq^{n(1+\epsilon)}).$$
\end{theorem}

In particular, there exists no $\delta > 0$ such that $\zeta_{\mathbf{H}({\Gri_S})}$ has a meromorphic extension to $\Re s > 1 - \delta$.

\section{Proof of Theorem \ref{thm1}}\label{subsec:proof.A2}

\subsection{Preliminaries}\label{prelim}
Recall that $\mathbf{H}$ is either an inner or an outer form. In the
latter case let $\chi$ denote the quadratic character of $k$
associated with the relevant extension $K/k$ describing the splitting
behaviour of a prime ideal $\mathfrak{p}$ of $k$ in $K$. In the former
case we denote by $\chi$ simply the trivial character.

If necessary, we momentarily enlarge the set $S$ to include the
finitely many finite places $v$ of $k$ (if any) where $K/k$ ramifies
or whose residue field cardinality is divisible by 2 or 3.

Let $v \not\in S$ be a place corresponding to a prime ideal
$\mathfrak{p}$ of $k$ with residue field cardinality $q=N
\mathfrak{p}$. As in \cite[(1.4)]{AKOV} we define
$$\varepsilon = \varepsilon_v  = \chi(\mathfrak{p}) \in \{-1,1\}.$$
If $\varepsilon = 1$, then
 $\mathbf{H}(\Gri_v) \cong {\rm SL}_3(\Gri_v)$; if $\varepsilon = -1$,
 then $\mathbf{H}(\Gri_v) \cong {\rm SU}_3(\Gri_v)$.  
 
 Let $\psi$ denote the  character defined by $$\psi_v = \psi(\mathfrak{p}) = \Big(\frac{-3}{q}\Big) = \begin{cases}   1, & q \equiv 1\, (\text{mod } 3),\\ - 1, & q \equiv - 1\, (\text{mod } 3). \end{cases}$$
  As in
 \cite[(1.10)]{AKOV} we define
$$ \iota = \iota(\varepsilon, q) = (q-\varepsilon, 3) = 2 + \epsilon_v \psi_v
 \in \{1, 3\}.$$  
Note that $\varepsilon$ and hence $\iota$ depend on $v$, not only on $q$. 

We write $\zeta_k^{S}(s) =\prod_{v\not\in S}\zeta_{k,v}(s)$ for the
zeta function $\zeta_k(s) = \prod_{v}\zeta_{k,v}(s)$ without the local
factors indexed by places $v\in S$, and similarly for the $L$-function
$L^S(s, \chi\psi)$.

In the function field case, the zeta- and $L$-functions under
consideration are periodic with respect to $s \mapsto s + \frac{2\pi
  i}{\log \qq}$, since its Dirichlet coefficients are indexed only by
powers of $\qq$.  We note that $\zeta_k$ has a simple poles at $s=0$
and $s=1$ (resp.\ $j + \frac{2\pi i}{\log \qq} n$, $n \in\Bbb{N}$, $j
\in \{0, 1\}$ in the function field case), but is holomorphic
otherwise. We need to understand the analytic behaviour of $L(.,
\chi\psi)$.
 
We first assume that $k$ is a number field. Then $\chi $ and $\psi$
are Hecke characters, and we observe that $\chi\psi$ is trivial if and
only we are in Case (A). Indeed, $\psi$ is trivial if and only if
$\zeta_3 \in k$. If $\psi$ is non-trivial, then $K/k$ is uniquely
determined by the condition $\chi = \psi$, and we see that for $K =
k(\zeta_3)$ the character of the extension $K/k$ equals~$\psi$.

We now assume that $k$ is a function field. Then $\chi$ is trivial if
and only if we are in case (A), namely $\mathbf{H}$ in an inner
form. Moreover, $\psi$ is the again the character associated with the
extension $k(\zeta_3)/k$, which is either trivial (if $\zeta_3 \in k$,
equivalently $\eta = 0$) or a quadratic constant field extension (if
$\zeta_3 \not\in k$, equivalently $\eta = 1$). From
\cite[Prop.~5.3.2]{friedjarden} we conclude that
\begin{equation}\label{Lpsi}
L(s, \chi\psi) = L\Big(s + \eta \frac{\pi}{\log \qq}, \chi\Big)
\end{equation}
 which is entire if and only if $\chi$ is trivial (since $K/k$ is not a constant field extension).



\subsection{Generic Euler factors}\label{subsubsec:non-arch}

 Our starting point for the proof of \Cref{thm1} is the Euler
 product~\eqref{equ:EP}.  
We use the explicit  formula \cite[Corollary D]{AKOV},  cf.\ \cite[Thm.~ C]{OPS/25} in the function field case, 
and consider  the product
 \begin{equation}\label{equ:euler.A2}
  \prod_{v \not \in S} \left(
   \zeta_{\mathbf{H}(\Bbb{F}_q)}(s) + \psi_{\varepsilon, q}(s) \right),
 \end{equation}
   where \begin{equation}\label{euler}
\begin{split}
  \zeta_{\mathbf{H}(\Bbb{F}_q)}(s) =&   1 + \frac{1}{(q^2 + \varepsilon q)^{s}} + \frac{q - 1 - \varepsilon}{(q^2 +  \varepsilon q + 1)^s} + \frac{q^2 - q - 1 + \varepsilon}{2(q^3 - \varepsilon)^s} + \frac{1}{q^{3s}} + \frac{q - 1 - \varepsilon }{(q^3 + \varepsilon q^2 + q)^s} \\
  &+ \frac{q^2 + \varepsilon q - 2 + 2 \iota(\varepsilon, q)^{2+s}}{3((q + \varepsilon)(q - \varepsilon)^2)^s}   + \frac{(q - \varepsilon)(q - 3 - \varepsilon) + 2\iota(\varepsilon, q)^{2 + s}}{6((q^2 + \varepsilon q+  1)(q + \varepsilon))^s} 
\end{split}
\end{equation}
 is the representation zeta function of the finite group of Lie type $\mathbf{H}(\Fq)$ and
\begin{equation}\label{euler1}
\begin{split}
 \psi_{\varepsilon, q}(s) = \Big(&\frac{(1 - q^{2 - 3s})(q-1)(q-\varepsilon)(2 + 2 q^{-s} + (q-2)(q + 1)^{-s} + q(q-1)^{-s})}{2 (q^2(q^2 + \varepsilon q + 1))^s}\\
 & + \frac{(1 - q^{2- 3s})(q - \varepsilon + \iota(\varepsilon, q)^{2 + s}(q + \varepsilon)(q - \varepsilon)^{-s} + \iota(\varepsilon, q)^2 (q - 1)(q^2 - 1)q^{-s})}{((q^3 - \varepsilon)(q + \varepsilon))^s}\\
 & + \frac{(q-1)(q - \varepsilon)^2 (q - 2 + 2 q^{2 - 2s} - q^{1 - 2s})}{6(q^3(q^2 + \varepsilon q + 1)(q + \varepsilon))^s} + \frac{(q-1)(q^2 - 1)q(1 - q^{-2s})}{2(q^3(q^3- \varepsilon))^s}\\
 & + \frac{(1 - q^{1 - 2s})(q^2 - 1)(q^2 + \varepsilon q + 1)}{3(q^3(q^2 - 1)(q - \varepsilon))^s}  + \frac{(q - 1)(q - \varepsilon)q ( 1 + q^{1 - 2s})}{(q^2 (q^3 - \varepsilon)(q  + \varepsilon))^s}\\
 & + \frac{(1 - q^{-2s})q^2 \iota(\varepsilon, q)^{2+s}}{(q(q^3 - \varepsilon)(q^2 - 1))^s} + \frac{(\varepsilon+ 1)\iota(\varepsilon, q)^{2+s} q^{2 - 2s}}{((q^3 - 1)(q^2 - 1) q)^s}\Big) \frac{1}{(1 - q^{1-2s})(1 - q^{2 - 3s})}.
 \end{split}
\end{equation}
We note that each of these Euler factors is a rational function in
finitely many numbers $n_i^{-s}$ (with $n_i$ depending on $q$,
$\varepsilon$, and $\iota$) and hence a meromorphic function. In $\Re
s > 1/2$ the only possible poles appear on the line $\Re s =
2/3$. Since $ \zeta_{\mathbf{H}({\Gri_S})}$ is a generating series of non-negative objects, it is non-zero on the segment $ s > 2/3$.

In order to make these terms resemble \eqref{poly}
more closely, we insert a Taylor expansion:
\begin{equation}\label{taylor}
\begin{split}
 & (q + \varepsilon)^ {-s} = q^{-s} \Big(1 - \frac{\varepsilon s}{ q} + O\Big(\frac{|s|^2}{q^2}\Big)\Big),\\
 &  (q^2 + \varepsilon q + 1)^{-s} = q^{-2s} \Big(1 - \frac{\varepsilon s}{ q} + O\Big(\frac{|s|^2}{q^2}\Big)\Big),\\
  &  (q^3 - \varepsilon)^{-s} = q^{-3s} \Big(1 +  O\Big(\frac{|s|^2}{q^2}\Big)\Big),\\
  & ((q + \varepsilon)(q - \varepsilon)^2)^{-s} =  q^{-3s} \Big(1 + \frac{\varepsilon s}{ q} + O\Big(\frac{|s|^2}{q^2}\Big)\Big),\\
   & ((q^2 + \varepsilon q + 1)(q \pm 1))^{-s} = q^{-3s}  \Big(1 - \frac{(  \varepsilon \pm 1) s }{q} + O\Big(\frac{|s|^2}{q^2}\Big)\Big),\\
     & ((q^3 - \varepsilon)(q^2 - 1))^{-s} = q^{-5s}\Big(1+ O\Big(\frac{|s|^2}{q^2}\Big)\Big), \\
    \end{split}
\end{equation}
\begin{equation*}
\begin{split}
   & ((q^2 - 1)(q- \varepsilon))^{-s} = q^{-3s} \Big(1 + \frac{ \varepsilon s}{q} + O\Big(\frac{|s|^2}{q^2}\Big)\Big),\\
  & ((q^3 - \varepsilon)(q + \varepsilon))^{-s} = q^{-4s}\Big(1 - \frac{\varepsilon s}{ q} + O\Big(\frac{|s|^2}{q^2}\Big)\Big).
  \end{split}
\end{equation*}


Plugging this into the previous equation, we obtain after a
straightforward computation
\begin{equation}\label{formula}
\begin{split}
\mathcal{E}(s, q) &:= (1 - q^{1-2s})(1 - q^{2 - 3s})\left(
\zeta_{\mathbf{H}(\Bbb{F}_q)}(s) + \psi_{\varepsilon, q}(s)\right) \\ &= 1 +
\iota(\varepsilon, q)^2(q^{3 - 5s} - q^{5 - 8s}) + 
F(s, q), 
\end{split}
\end{equation}
say, with
\begin{equation}\label{boundF}
F(s, q) \ll_s q^{4 -8 \Re s}  + q^{2 - 5\Re s} + q^{- 2\Re s}  \ll q^{4 -8 \Re s}   + q^{- 2\Re s}
\end{equation}
and
\begin{equation*}
\frac{d}{ds} F(s, q)  = \int_{|z - s| = (\log q)^{-1}}  \frac{F(z, q)}{z-s} \frac{dz}{2\pi i}  \ll_s (q^{4 -8 \Re s}  + q^{- 2\Re s})\log q .
\end{equation*}

\begin{rem}
The summands in~\eqref{euler1} reflect the organization of the
relevant characters according to invariants called \emph{shadows} in
\cite{AKOV}. Our analysis shows that no single shadow or summand
in~\eqref{euler1} suffices to explain the asymptotic properties
established in~\Cref{thm1}. The relatively simple shape of
\eqref{formula} is therefore quite remarkable: both \eqref{euler} and
\eqref{euler1} contribute additional terms of the form $q^{4-6s}$, but
they cancel. If the did not, we could only continue to $\Re s > 2/3$,
cf.\ also Remark \ref{rem21}.
\end{rem}


In the following three sections we show that the infinite product
\eqref{equ:euler.A2} has meromorphic continuation to $\Re s > 5/8$ and
a natural boundary at $\Re s = 5/8$. In \Cref{subsubsec:except} we
show that these properties remain true after adding the finitely many
missing Euler factors.
 
 \subsection{The analytic continuation}
 We now start with the proof of the analytic continuation to $\Re s > 5/8$. 
For $n, m \in \Bbb{Z}$ with $n \equiv m \,  (\text{mod }2)$ 
we define the polynomials
\begin{align*}
  P_+(x; n, m) &= \begin{cases} (1 - x)^n, & |m| \leq n,\\ (1 + x)^{-(n+m)/2} (1 - x)^{(n-m)/2},  &|n| \leq -m,\\(1+x)^{-n}, & |m| \leq -n,\\ (1 - x)^{(n+m)/2} (1+ x)^{-(n-m)/2},  &|n| \leq m,\end{cases}\\
  P_-(x; n, m) &= \begin{cases}  (1 - x)^{(n+m)/2} (1 + x)^{(n-m)/2} , & |m| \geq n,\\ ( 1 + x)^{-m}, & |n| \leq -m\\(1+ x)^{-(n+m)/2} (1 - x)^{-(n-m)/2} , & |m| \leq -n,\\ ( 1 - x)^{m}, & |n| \leq m.\end{cases}
\end{align*}

Clearly we have $P_+(x; n, m) \equiv P_-(x;n, m)$ (mod $2$) and also
\begin{equation}\label{linear}
P_+(x; n, m) =  1 - nx + \ldots,  \quad P_-(x;n, m) = 1 - m x + \ldots. 
\end{equation}
By direct comparison of Euler products we obtain
\begin{equation}\label{zetaL}
\begin{split}
&\prod_{ \substack{v \not\in S \\ \varepsilon_v\psi_v = 1}} P_+(q^{-s}; n, m)\prod_{ \substack{v \not\in S \\ \varepsilon_v\psi_v= -1}}P_-(q^{-s}; n, m)   \\ &=\zeta_k^S(s)^{-\frac{n+m}{2}} L^{S}(s, \chi \psi)^{-\frac{n-m}{2}}  \cdot \begin{cases} 1, & |m| \leq n,\\ \zeta_k^S(2s)^{\frac{n+m}{2}},  &|n| \leq -m,\\\zeta_k^S(2s)^{n}, & |m| \leq -n,\\ \zeta_k^S(2s)^{\frac{n-m}{2}}. &|n| \leq m.\end{cases}
\end{split}
\end{equation}


We now return to the right hand side of \eqref{formula} and modify the ideas of \cite[Lemma 5.5]{dSW}. We put $X = q$, $Y = q^{-s}$ and consider the polynomials 
$$W_{+, 0}(X, Y) = 1 + 9 X^3 Y^5 - 9 X^5 Y^8, \quad W_{-, 0}(X, Y) = 1
+ X^3 Y^5 - X^5 Y^8.$$ We   proceed recursively as follows: given
two polynomials $$W_{+, j} = \sum_{m, n} \alpha^{(j)}_{m, n} X^n Y^m,
\quad W_{-, j} = \sum_{m, n} \beta^{(j)}_{m, n} X^n Y^m\in \Bbb{Z}[X,
  Y]$$ with $\alpha^{(j)}_{m, n} \equiv \beta^{(j)}_{m, n}$ (mod $2$)
and $\alpha^{(j)}_{0, 0} = \beta^{(j)}_{0, 0} = 1$, we order the
monomials lexicographically by their exponents $(m,n)$ and pick the
smallest index $(m_j, n_j) > (0, 0)$ such that $\alpha^{(j)}_{m_j,
  n_j}$ or $\beta^{(j)}_{m_j, n_j}$ are non-zero. We define
\begin{align*}
  W_{\pm, j+1}(X, Y) &= W_{\pm, j}(X, Y) P_{\pm}(X^{n_j} Y^{m_j};
  \alpha^{(j)}_{m_j, n_j}, \beta^{(j)}_{m_j, n_j}).
\end{align*}
Then the polynomials $W_{\pm, j+1} $ have again constant term 1 and
are congruent modulo $2$. Moreover, the smallest nontrivial monomial
$X^{n_j} Y^{m_j}$ is cleared by \eqref{linear} and no smaller
monomials are inferred. Finally, inductively we see easily that only
monomials of the form
\begin{equation}\label{exponents}
 (X^3 Y^5)^u    (X^5 Y^8)^v = X^{3u + 5v} Y^{5u + 8v}
 \end{equation}
for $u, v \in \Bbb{N}_0$  can occur in $W_{{\pm}, j}$. 


For illustration we carry out the first two steps of this procedure. We have $(m_0, n_0) = (5, 3)$, $\alpha^{(0)}_{5, 3} = 9$, $\beta^{(0)}_{5, 3} = 1$ and 
\begin{equation}\label{step1a}
P_{\pm}(X^3 Y^5, 9, 1) = \begin{cases} (1 - X^3 Y^5)^9, & \pm = +,\\ (1 - X^3 Y^5)^5(1 - X^3 Y^5)^4, & \pm = -,  \end{cases}\end{equation}
and so
\begin{equation}\label{step1b}
W_{\pm, 1} = \begin{cases}  1 - 9 X^5 Y^8 - 45 X^6 Y^{10} + \ldots  + 9 X^{32} Y^{53}, & \pm = +,\\ 1 -  X^5 Y^8 -  5 X^6 Y^{10} + \ldots  +  X^{32} Y^{53}, & \pm = -. 
\end{cases}
\end{equation}
Next, we have $(m_1, n_1) = (8, 5)$ and $\alpha^{(1)}_{8, 5} = -9$, $\beta^{(1)}_{8, 3} = -1$ and 
$$P_{\pm}(X^5 Y^8, -9, -1) = \begin{cases} (1 + X^5 Y^8)^9, & \pm = +,\\ (1 + X^5 Y^8)^5(1+ X^5 Y^8)^4, & \pm = -,  \end{cases}$$
getting
$$W_{\pm, 2} = \begin{cases}  1  - 45 X^6 Y^{10} + \ldots  + 9 X^{77} Y^{125}, & \pm = +,\\ 1  -  5 X^6 Y^{10} + \ldots  +  X^{77} Y^{125}, & \pm = -. 
\end{cases}$$


For $v \not\in S$ and $i \in \Bbb{N}$ let us define
$$\mathcal{P}_{i, v} := \prod_{j = 0}^{i-1} P_{\varepsilon_v \psi_v }(q^{n_j + sm_j}; \alpha^{(j)}_{m_j, n_j}, \beta^{(j)}_{m_j, n_j}). $$
Clearly, we have $$W_{\pm, 0}(X, Y) = W_{\pm, i}(X, Y) \prod_{j = 0}^{i-1}   P_{\pm}(X^{n_j} Y^{m_j}; \alpha^{(j)}_{m_j, n_j}, \beta^{(j)}_{m_j, n_j})^{-1} $$
for any $i \in \Bbb{N}$, and by \eqref{zetaL} and \eqref{exponents} we know that 
\begin{equation*}
  \prod_{\pm} \prod_{\substack{v \not \in S \\ \varepsilon_v \psi_v = \pm 1} } \mathcal{P}_{i, v}^{-1}
  \end{equation*}
is a finite product   of positive or negative integral powers of
\begin{equation}\label{zetaL1}
\zeta_k^S( (5u+8v) s - (3u + 5v) )  \quad \text{and}\quad  L^{S}( (5u+8v) s - (3u + 5v) , \chi\psi) 
\end{equation}  for certain $u, v \in \Bbb{N}_0$, $(u, v) \not= (0, 0)$, in particular meromorphic in $s$.  

Returning to \eqref{formula}, we can write
\begin{equation}\label{corr}
\begin{split}
\mathcal{E}(s, q) = W_{\varepsilon_v\psi_v, 0} (q, q^{-s})   + F(s, q) = \big(W_{\varepsilon_v\psi_v, i} (q, q^{-s})  + F(s, q) \mathcal{P}_{i, v}\big)\mathcal{P}_{i, v}^{-1} 
\end{split}
\end{equation}
for any $i \in \Bbb{N}$. Here 
$$W_{\varepsilon_v\psi_v, i} (q, q^{-s})  = 1 + q^{(3u_0 + 5v_0) - (5u_0 + 8v_0)s} + \ldots$$
where $(u_0, v_0)$ can be chosen as large as we wish by choosing $i$ sufficiently large. Moreover, by \eqref{boundF} we see that 
$F(s, q)\mathcal{P}_{i, v}$ 
is bounded by
$$q^{(4 + 3u + 5v) - (8 + 5u+8v) \Re s} + q^{( 3u + 5v) - (8 + 5u+8v + 2) \Re s}$$
for certain positive integers  $u, v$. We conclude that
\begin{equation}\label{euler-conv}
\prod_{\pm} \prod_{\substack{v\not\in S \\\varepsilon_v\psi_v = \pm 1}}\big(W_{\varepsilon_v\psi_v, i} (q, q^{-s})  + F(s, q)\mathcal{P}_{i, v}
\big)
\end{equation}
is absolutely convergent in  $\Re s \geq 5/8 + \delta$ for any $\delta > 0$. 
In this way we obtain a meromorphic continuation of $$\zeta_{\mathbf{H}({\Gri_S})}(s) =\zeta_k^S( 2s-1)\zeta_k^S( 3s-2)  \prod_{v\not\in S} \mathcal{E}(s, q)$$ to the half plane $\Re s > 5/8$. \\

Carrying out only the first step of the inductive procedure described above, we see from \eqref{step1a} and \eqref{step1b} that
\begin{equation}\label{step1}
 \zeta_{\mathbf{H}({\Gri_S})}(s) = \zeta_k^{S}( 2s-1)\zeta_k^{S}( 3s -2) \zeta_k^{S}(5s - 3)^5 L^{S}(5s - 3, \chi\psi)^{4 } H (s)
 \end{equation}
where $H $ is an absolutely convergent Euler product in $\Re s > 3/4$ and hence in particular holomorphic and non-vanishing. 

Thus in $\Re s > 3/4$ the  function $ \zeta_{\mathbf{H}({\Gri}_S)}$ has the polar behaviour described in Theorem \ref{thm1}.

\subsection{The natural boundary -- number field case}
Next we show that the line $\Re s = 5/8$ is a natural boundary. To
this end, we show that every given point $s_0 = 5/8 + it$ is a limit
point of zeros. For clarity we assume in this subsection that $k$ is a
number field and explain the necessary modifications in the function
field case in the next subsection. In both cases we show, roughly
speaking, that sufficiently many Euler factors have zeros in a small
neighbourhood of $s_0$. A difficulty is, however, to make sure that
these zeros are not cancelled by poles for the zeta- and $L$-functions
that arise in the course of the analytic continuation. Here the
argument diverges in the number and function field case. In the former
we use bounds for the total number of zeros of zeta functions over
number fields. In the latter we use the Riemann hypothesis (which is
known over function fields) along with Kronecker's simultaneous
approximation theorem to compensate for the fact that we have only
very few different norms at our disposal (namely those that are powers
of $\mathfrak{q}$).


We proceed to show that $\Re s= 5/8$ is a natural boundary in the
number field case.  Let $c_0 = c_0(t) = 9(\pi + |t| \log 2)/\log 9$.
Fix some small $\delta > 0$ and consider the
rectangle $$\mathcal{R}_{\delta} = \{ s \in \Bbb{C} : |\Re (s - s_0
-\delta)| \leq \delta/2, |\Im (s - s_0)| \leq c_0 \delta \}$$ with a
typical Euler factor $\mathcal{E}(s, q)$ as in \eqref{formula}
corresponding to a place $v$ with $\varepsilon_v\psi_v = 1$. We recall
that in this case
\begin{equation}\label{eulerrecall}
\mathcal{E}(s, q) = 1 + 9q^{3 - 5s} - 9q^{5 - 8s} + O(q^{4 - 8\Re s } + q^{2- \Re s}).
\end{equation}
In what follows we always assume that $q$ is sufficiently large in
terms of $t$ and the fixed field extension $K/k$. We will choose later
$q$ to be a function of $\delta$ and let $\delta$ tend to zero.

For $n \in \Bbb{Z}$ we are looking for zeros of $\mathcal{E}(s, q)$ in
a neighbourhood of the points
$$s=  \frac{5}{8} + \frac{i(1 + 2n)\pi }{\log q}. $$
A good approximation can be found by putting $V = q^{-1/8}$, $U = q^{5/8 - s}$ and writing down the Puiseux series in $V$ of the equation $1 + 9U^5V - 9U^8 = 0$ near $U = -3^{-1/4}$. One checks that
$$1 + 9q^{3 - 5s} - 9q^{5 - 8s}|_{s = \frac{5}{8} +  \frac{1}{\log q}  (i(1 + 2n)\pi + \frac{\log 9}{8})} \ll q^{-1/8}$$
and hence there must be constants $c_j$ such that for 
\begin{equation}\label{defsnq}
\begin{split}
s_{n, q} := \frac{5}{8} +  \frac{1}{\log q} \Big(&i(1 + 2n)\pi + \frac{\log 9}{8} +  \sum_{j=1}^7 c_j q^{-j/8}\Big)
\end{split}
\end{equation}
we have 
$$1 + 9q^{3-5s_{n, q}} - 9q^{5 - 8s_{n, q}} \ll q^{-1}. $$
(While not relevant for the following discussion, the constants are
\begin{displaymath}
\begin{split}
& c_1 =  \frac{3^{3/4}}{8} , \quad c_2 =  - \frac{3\sqrt{3}}{64} , \quad c_3 = - \frac{21 \cdot 3^{1/4}}{1024}, \quad c_4 = \frac{27}{512} , \quad c_5 =  - \frac{9639\cdot  3^{3/4}}{1310720},\\
& c_6 = - \frac{2079 \sqrt{3}}{131072}, \quad c_7 =  \frac{5942079 \cdot 3^{1/4}}{234881024}. \Big)
\end{split}
\end{displaymath} 
Note that $\Re s_{n, q} > 5/8$. We assume $n \ll \log q$, so that $s_{n, q} \ll 1$. 
By  \eqref{eulerrecall} we obtain
$$\mathcal{E}(s_{n, q}, q) 
\ll q^{-1}. $$
On the other hand, for $|s  - s_{n, q}| \ll q^{-1}$ we compute
$$\frac{d}{ds} \mathcal{E}(s, q) = (8 q^{5 - 8s} - 5q^{3 - 5s})\log q
+ \frac{d}{ds} F(s, q) = 8 + O(q^{-1/8} \log q).$$ Thus for $q$
sufficiently large and $n \ll \log q$ we find a point $s_{n,
  q}^{\ast}$ (for instance by Newton's method) with 
  \begin{equation}\label{approx}
     s_{n, q}^{\ast} -
s_{n, q} \ll q^{-1}
\end{equation}
 such that $\mathcal{E}(s_{n, q}^{\ast}, q) = 0$. Note 
that for $q \leq q' < 2q$ we have
\begin{equation}\label{distinct}
\begin{split}
|s^{\ast}_{n, q} - s^{\ast}_{n, q'} | &\geq |\Im (s^{\ast}_{n, q} - s^{\ast}_{n, q'} )|  = (1 + 2n)\pi \Big( \frac{1}{\log q} - \frac{1}{\log q'}\Big) + O\Big(\frac{1}{q}\Big) \\
&\asymp \frac{|1 + 2n|(q' - q)}{q (\log q)^2} + O\Big(\frac{1}{q}\Big)  > 0
\end{split}
\end{equation}
provided that $q' - q \geq (\log q)^3$. (Note that because of the error term $F(s, q)$ we cannot use algebraic arguments such as  \cite[p.\ 130 or p.\ 133]{dSW}.) 

We now define $q_0$ and $n$ by
\begin{equation}\label{defn}
\delta = \frac{\log 9}{8\log q_0}, \quad n = \Big[ \frac{t \log q_0}{2\pi}\Big]
\end{equation}
and restrict to numbers $q \in \mathcal{I} = [q_0, 2q_0]$. 
Then for $q \in \mathcal{I}$ we have 
$$\Big|\frac{(1+2n)\pi}{\log q} - t\Big| \leq \frac{\pi + |t| \log 2}{\log q_0} = \frac{8}{9} c_0 \delta$$ 
and we conclude from  
\eqref{approx}   that for sufficiently large $q$ we have 
\begin{displaymath}
\begin{split}
&|\Im s_{n, q}^{\ast} - t | \leq c_0 \delta, \quad  \Big|\Re s_{n, q}^{\ast} - \frac{5}{8}  - \delta \Big| \ll \delta^2 \leq \frac{1}{2}\delta \end{split}
 \end{displaymath}
and so $$s_{n, q}^{\ast} \in \mathcal{R}_{\delta}.$$ By
\eqref{distinct} all $s_{n, q}^{\ast}$ with $q \in \mathcal{I}$ are
pairwise distinct, provided the values of $q$ are at least $(\log
q_0)^3$-spaced. We can force this by restricting to places $v$ such
that $q \equiv 1$ (mod $\lceil (\log q_0)^3 \rceil$). We recall in
addition the condition $\epsilon_v\psi_v = 1$, which is another set
of (fixed) congruences modulo $3\cdot {\rm Nr}(\Delta)$ where $\Delta
\in k$ is the discriminant of the extension $K/k$. The number of such
places can be evaluated by a number field version of the
Siegel-Walfisz theorem \cite{Mi} with the required uniformity in the
modulus: we find
$$\asymp_{K/k} \frac{q_0}{\log q_0} \cdot \frac{1}{\log^{3[k: \Bbb{Q}]} q_0} \gg
q_0^{0.9} \geq \exp(\tfrac{1}{10} \delta^{-1})$$   places $v$ with $q \in \mathcal{I}$  such that $\mathcal{E}(s, q)$ has a zero in $\mathcal{R}_{\delta}.$ 

We now return to \eqref{corr} and write 
$$\mathcal{E}(s, q) = \big(\mathcal{E}(s, q)\mathcal{P}_{i, v}\big)  \cdot  \mathcal{P}_{i, v}^{-1} $$
for  $i\in \Bbb{N}$. By the discussion around \eqref{euler-conv}, the global Euler product $\prod_{v\not \in S}\mathcal{E}(s, q)\mathcal{P}_{i, v}$ is absolutely convergent in $\mathcal{R}_{\delta}$ upon choosing $i$ sufficiently large in terms of $\delta$, and we know that  at least $\exp(\frac{1}{10}\delta^{-1})$ Euler factors corresponding to places $v$ with $q\in \mathcal{I}$  have a zero  $\rho \in \mathcal{R}_{\delta}$, and all these zeros are distinct. 

On the other hand, we can express the global Euler product $\prod_{v\not \in S}   \mathcal{P}_{i, v}^{-1}$ as a finite product   of positive or negative integral powers of zeta- and $L$-functions of the form \eqref{zetaL1}.   Suppose that some of the above zeros $\rho \in \mathcal{R}_{\delta}$ coincides with  a zero
of some $\zeta$- or $L$-factors in \eqref{zetaL1} (which may appear in the denominator with some multiplicity and hence cancel $\rho$). Then
$$\frac{5}{8} + \delta + O(\delta^2) = \Re(\rho) \leq \frac{1 + 3u +
  5v}{5u + 8v} \leq \frac{1}{5u + 8v} + \frac{5}{8}$$\\ and $\Im \rho
\ll 1$, so that $$(5u + 8v)\rho - (3u + 5v)$$ is a zero of $\zeta_k$
or $L(., \chi\psi)$ with imaginary part $\ll (5u + 8v) \ll
1/\delta$. There are at most $O(\delta^{-1} |\log \delta|)$ such zeros
(cf.\ \cite[Thm~5.8]{IK}), and the number of choices for the pair $(u,
v)$ is $O(\delta^{-2})$, so that in total the meromorphic continuation
of the global Euler product $\prod_{v\not \in S} \mathcal{P}_{i,
  v}^{-1}$ can have at most $O(\delta^{-3} |\log \delta|)$ poles in
$\mathcal{R}_{\delta}$. Choosing $\delta$ sufficiently small, this
cannot compensate the zeros found above.

We conclude that we find a sequence of zeros of the product \eqref{equ:euler.A2} 
converging to our given point $s_0$.



 
 \subsection{The natural boundary -- function field case} The main difference in the function field case is that we have only a very sparse set of distinct values $q$ available, namely only numbers of the form $\qq^m$, $m\in \Bbb{N}$.  On the other hand, the Riemann hypothesis is known, which we will use in the subsequent analysis.
 
 As before, let us fix some $s_0 = 5/8 + it$ and choose some place $v$ with $q = \qq^m$ and $\varepsilon_v\psi_v = 1$. At least if $m$ is even and sufficiently large, this is always possible: for even $m$ we have $\psi_v = 1$ and by a basic form of the Chebotarev density theorem (e.g.\ \cite[Proposition 7.4.8]{friedjarden}  with $n=1$, $m=2$) places $v$ with $\epsilon_v=1$ exist in abundance. Let us choose one such $v$ for each such $q = \qq^{2m}$. 
  Choose $s_{n, q}$ as in \eqref{defsnq}, $s^{\ast}_{n, q}$ as in \eqref{approx}  and $n$ as in \eqref{defn}, so that $\mathcal{E}(s^{\ast}_{n, q}, q) = 0$ and 
\begin{equation}\label{defeta}
|\Im s_{n, q}^{\ast} - t |\ll \frac{1}{\log q}, \quad \Re s^{\ast}_{n, q} - \frac{5}{8} =:\eta_{q} = \frac{1}{\log q}\Big(\frac{\log 9}{8} + O(q^{-1/4})\Big).
\end{equation}
 For $q = {\qq}^{2m}$ this gives a sequences of distinct zeros tending to $s_0$. We need to show that it contains a subsequence that is not cancelled by possible zeros of the zeta- and $L$-factors at arguments of the form $(5u + 8v) s- (3u + 5v)$ with $u, v \in \Bbb{N}_0$. The Riemann hypothesis (see e.g.\ \cite[Thm~5.5.1]{friedjarden}) states that the zeros of  $\zeta_K = \zeta_k L(., \chi)$ and hence of both $\zeta_k$ and $L(., \chi\psi)$ (cf.\ \eqref{Lpsi}) are on the line $\Re s = 1/2$. Therefore, the real parts of the zeros in question are at 
 $$\frac{\frac{1}{2} + 3u + 5v}{5u+8v}.$$
 Elementary algebra shows that this is larger than $5/8$ only for $u < 4$, and so it can only coincide with $5/8 + \eta_q$ if $(u, v)$ equals
 $$\Big(0, \frac{1}{16\eta_q}\Big), \quad\Big(1, \frac{3}{64\eta_q} - \frac{5}{8}\Big), \quad \Big(2, \frac{1}{32\eta_q} - \frac{5}{4}\Big), \quad \text{or}\quad \Big(3, \frac{1}{64\eta_q} - \frac{15}{8}\Big).$$
 In order to derive at a contradiction, we show that for an infinite subsequence of $q = \qq^{2m}$ the second entry is not an integer, say its fractional part $\| . \|$ is in $(1/4, 3/4)$. 
  To this end let $\alpha:= \frac{\log \qq}{2\log 9}$ and consider the linear polynomials 
  $$\ell_1(m) = 2\alpha m + \frac{1}{2}, \quad \ell_2(m) = \frac{3}{2}\alpha m   - \frac{1}{8}, \quad \ell_3(m) =  \alpha m - \frac{3}{4}, \quad \ell_4(m) = \frac{1}{2}\alpha m - \frac{11}{8},$$
which by the definition \eqref{defeta} of $\eta_q$ describe  the possible values of $v + 1/2$   up to an error of size $O(\qq^{-2m/4})$. 
Since $\alpha$ is irrational (since $3$ and $\mathfrak{q}$ are coprime), it follows from Kronecker's approximation theorem (see e.g.\ \cite[Thm.~442]{HW}) that there are infinitely many $m$ such that simultaneously $\| \ell_j(m) \| < 1/4$ for $1 \leq j \leq 4$, which by the above discussion provides the desired sequence of zeros tending to $s_0$. 

\subsection{The remaining Euler factors}\label{subsubsec:except}
Recall that at the beginning of the proof we have possibly enlarged
the set $S$ by finitely many places. We now deal with the local
factors at these places and the infinite Euler factor.  By
\cite[Thm.~1.1]{JZ} and \cite[Thm.~B]{AKOV1} the associated
local representation zeta functions $\zeta_{\mathbf{H}(\Gri_v)}(s)$
are rational functions (hence meromorphic on the whole complex plane)
with their right-most poles at $\Re s = 2/3$. As generating function
of non-negative numbers, they are non-vanishing on the segment $ s >
2/3$, so they do not change the polar behaviour there, nor do they
change the natural boundary.

We finally discuss the Archimedean Euler factors of \eqref{equ:EP},
which only occur in the number field case.  They are given
(cf.\ \cite[p.\ 359]{KMT}) by powers of
$$\zeta_{{\rm SL}_3(\Bbb{C})}(s) = \sum_{m, n} \frac{1}{m^s
  n^s(m+n)^{s}}.$$ It suffices to obtain meromorphic continuation past
the region $\Re s > 2/3$ of absolute convergence. This is
straightforward (and well-known): by the integral formula for the beta function
\cite[3.196.2 with $u = 0$]{GR} we have, initially for $\Re s > 2/3$,
the absolutely convergent expression
\begin{displaymath}
\begin{split}
  \zeta_{{\rm SL}_3(\Bbb{C})}(s) = \sum_{n, m} \frac{1}{m^{2s}n^s (1 + n/m)^s} &= \sum_{m, n} \frac{1}{m^{2s}n^s} \int_{(1/3)} \frac{\Gamma(s-t) \Gamma(t)}{\Gamma(s)} \Big( \frac{n}{m}\Big)^{-t} \frac{dt}{2\pi i}\\
  & =  \int_{(1/3)} \frac{\Gamma(s-t) \Gamma(t)}{\Gamma(s)} \zeta(s+t) \zeta(2s - t)\frac{dt}{2\pi i}.
\end{split}
\end{displaymath}
Let us assume $2/3 < \Re s < 1$ and fix $\varepsilon > 0$ sufficiently small. Then shifting the contour to $\Re t = \varepsilon$, we obtain
\begin{displaymath}
\begin{split}
  \zeta_{{\rm SL}_3(\Bbb{C})}(s) =  \int_{(\varepsilon)} \frac{\Gamma(s-t) \Gamma(t)}{\Gamma(s)} \zeta(s+t) \zeta(2s - t)\frac{dt}{2\pi i} + \frac{\Gamma(1-s)\Gamma(2s - 1)}{\Gamma(s)} \zeta(3s - 1).
\end{split}
\end{displaymath}
This is meromorphic in $1/2 < \Re s < 1$ with a pole only at $s = 2/3$.

\begin{rem}\label{rem21}  One can form the global Euler product
 $$\prod_{v \not\in S} \zeta_{\textbf{H}(\Bbb{F}_q)}(s),$$
 and by a similar Taylor expansion obtain
 $$ \zeta_{\textbf{H}(\Bbb{F}_q)}(s)  = 1 + q^{2 - 3s} + q^{1 - 2s} + F^{\ast}(s, q)$$
 with $F^{\ast}(s, q) \ll q^{-2\Re s}$. Arguing similarly as before, we can  extend the global zeta function  meromorphically only to $\Re s > 2/3$. It is interesting to note that this half plane of continuation is smaller than the half plane of  continuation for  $\zeta_{\textbf{{\rm \bf{H}}}(\Gri_S)}$. In other words, as observed before, there is some non-trivial cancellation in the local factors $ \zeta_{\textbf{H}(\Bbb{F}_q)}(s)$ and $ \psi_{\varepsilon, q}(s) $. 
 \end{rem}
 
\section{Proof of \Cref{thm2}}
 This is a standard application of Mellin inversion. The technical
 difficulty lies in the fact that it is not clear (and quite possibly
 not true) that in the number field case $\zeta_{\textbf{{\rm \bf{H}}}(\Gri_S)}$ has
 polynomial growth on vertical lines beyond the region of absolute
 convergence, cf.\ the remark after \Cref{thm2}. This is irrelevant in the function field case since all zeta functions are $\frac{2\pi i}{\log \qq}\Bbb{Z}$-periodic. 
 
 We first give the argument in the number field case. To obtain subexponential growth in the region $\Re s \geq 3/4$, we
replace all parentheses on the right hand sides of \eqref{taylor} by
$$ 1 + O\Big(\min\Big(1, \frac{|s|}{q}\Big)\Big) = 1 +
O\Big(\frac{|s|^{99/100}}{q^{99/100}}\Big). $$ Then by the same
computation we obtain
$$F(s, q) \ll |s|^{99/100} \big(q^{1 - 3\Re s} +  q^{-2\Re s}\big)$$
in $\Re s \geq 3/4$ where now the implied constant is absolute. For fixed $0 < \epsilon < 1/100$, we can then  write \eqref{step1} as
$$ \zeta_{\mathbf{H}({\Gri}_S)}(s) = \zeta_k^{S}( 2s-1)\zeta_k^{S}( 3s -2) \zeta_k^{S}(5s - 3)^5 L^{(S)}(5s - 3, \chi \psi)^{4 } H (s)
$$
where
$$H(s) = H_0(s)\prod_{v \not \in S} \Big(1 + \frac{O(1)}{q^{1 + \epsilon}} + \frac{O(|s|^{99/100})}{q^{5/4}}\Big) \ll \exp(|s|^{99/100})$$
in $\Re s \geq 3/4 + \epsilon$ where $H_0$ contains potentially finitely many Euler factors that we discard at the beginning of Section \ref{prelim}. 
By Mellin inversion (cf.\ \cite[(4.107)]{IK}) 
we now obtain
$$\sum_{n } r_n(\textbf{{\rm \bf{H}}}(\Gri_S)) e^{-n/x} = \int_{(2)}  \zeta_{\textbf{H}(\Gri_S)}(s) \Gamma(s) x^s \frac{ds}{2\pi i}.$$
We shift the contour to $\Re s = 3/4 + \epsilon$. The residues at $s = 1$ and $s = 4/5$ yield the two main terms, and the remaining integral is rapidly convergent on the line $\Re s = 3/4 + \epsilon$ in view of the estimate $\Gamma(s) \ll (1 + |s|)^{1/2}\exp(-\pi  |s|/2)$. This completes the proof. 
\medskip

If $k$ is a function field with constant field $\Bbb{F}_{\qq}$, then we have the Mellin formula
$$  r_{\qq^n}(\textbf{{\rm \bf{H}}}(\Gri_S))  = \int_{2 - \frac{\pi i}{\log \qq}}^{2 + \frac{\pi i}{\log \qq}}  \zeta_{\textbf{H}(\Gri_S)}(s)   (\qq^n)^s (\log \qq) \frac{ds}{2\pi i}$$
which in this case is nothing but Cauchy's integral formula for the power series $\sum_n   r_{\qq^n}(\textbf{{\rm \bf{H}}}(\Gri_S)) x^n$ with $x = \qq^{-s}$. We can now shift the contour to $\Re s \geq 3/4 + \epsilon$ in the same way as before, noting that the contribution of the horizontal lines cancel by periodicity.

\section{Groups of type $\mathsf{A}_1$ -- proof of \Cref{prop3}}

This follows from the local formulae \cite[(A.5)]{AKOV}, see also
\cite[Thm.~7.5]{JZ}, which read
\begin{displaymath}
\begin{split}
\zeta_{\textbf{H}(\Gri_v)}(s) = &1 + q^{-s} + \frac{(q-3)}{2(q+1)^s} + \frac{(q-1)^{1-s}}{2} + 2 \Big( \frac{q+1}{2}\Big)^{-s} + 2\Big(\frac{q-1}{2}\Big)^{-s}\\
& + \frac{2^{2+s}q(q^2 - 1)^{-s} + \frac{1}{2}(q-1)(q- \varepsilon)(q^2 - \varepsilon q)^{-s} + \frac{1}{2}(q-1)(q- \varepsilon) (q^2 + \epsilon q)^{-s}}{1 - q^{1-s}}
\end{split}
\end{displaymath}
for places $v$ of norm $q \geq 5$ where $\textbf{H}$ is ${\rm SL}_2$ or ${\rm SU}_2$ and accordingly $\varepsilon = 1$ resp.\ $-1$. 
By a similar Taylor argument as before we obtain
$$\zeta_{\textbf{H}(\Gri_v)}(s) = (1 - q^{1-s})^{-1} \big(1 + O_s(q^{1-2\Re s} + q^{-\Re s})\big)$$
and the meromorphic continuation to $\Re s > 1$ follows easily. The asymptotic formula follows as in the previous proof.

To analyze the Euler product further as we approach $\Re s = 1$, we continue with   Taylor expansions, and for $\Re s > 1/2$ we   write 
 $$1 - \frac{\varepsilon}{q^s} + \frac{2^{2+2s} - 1}{q^{2s-1}}   =  \Big(1 + \frac{\epsilon}{q^s} \Big)^{-1}\Big(1- \frac{1}{q^{2s-1}}\Big)^{1-2^{2+s}  } \Big(1   + O_s\big( q^{2 - 4\Re s} +q^{-2\Re s}\big)\Big).$$
If we write as before $\chi(\mathfrak{p}) = \varepsilon_v$ for a finite place $v = \mathfrak{p}$, we obtain
$$\zeta_{\mathbf{H}({\Gri_S})}(s) = \zeta(s)^{r}
\frac{\zeta_k^S(s-1)}{L^S(s, \chi)} \zeta^S_k(2s-1)^{ 2^{2+s}-1}
H(s) $$ (with $r=0$ in the function field case and $r = [k:\Bbb{Q}]$
in the number field case) for $ \Re s > 1$ where $H$ is holomorphic
and nonzero in $\Re s > 3/4$. As $\zeta_k$ has a pole at $s=1$, the
factor $\zeta_k(2s-1)^{ 2^{2+s}-1}$ and hence
$$\zeta^S_k(2s-1)^{ 2^{2+s}-1}$$
has a branch cut singularity at $s = 1$. More precisely, as $\zeta_k(s) = c(s-1)^{-1} + O(1)$ as $s \rightarrow 1$ 
for a constant 
$c  \not=0$, we have 
$$\zeta_k(2s-1)^{ 2^{2+s}-1} = \frac{c^7}{128(s-1)^7} 
+ \frac{c^7 \log 2}{16(s-1)^6} \log\Big( \frac{\tfrac{1}{2}c}{(s-1)}\Big) + O(|s-1|^6).$$
The factor $\zeta(s)^r/ L^S(s, \chi)$ may change the exponent $7$, but does not affect the type of the brunch cut singularity. 

\begin{acknowledgements}
We thank J\"urgen Kl\"uners, Uri Onn and Sun Woo Park  for very helpful conversations and discussions. Both authors were funded by the Deutsche Forschungsgemeinschaft
(DFG, German Research Foundation) -- Project-ID 491392403 -- TRR~358.
\end{acknowledgements}

\end{document}